\documentclass{amsart}
\title{On the indices of curves over local fields}
\author{Pete L. Clark}
\email{pete@math.uga.edu}

\DeclareFontEncoding{OT2}{}{} 
  \newcommand{\textcyr}[1]{%
    {\fontencoding{OT2}\fontfamily{wncyr}\fontseries{m}\fontshape{n}%
     \selectfont #1}}
\newcommand{\Sha}{{\mbox{\textcyr{Sh}}}}

\begin{document}
\newtheorem{lemma}{Lemma}
\newtheorem{prop}[lemma]{Proposition}
\newtheorem{cor}[lemma]{Corollary}
\newtheorem{thm}[lemma]{Theorem}
\newtheorem{ques}[lemma]{Question}
\newtheorem{quest}[lemma]{Question}
\newtheorem{conj}[lemma]{Conjecture}
\newtheorem{fact}[lemma]{Fact}
\newtheorem*{mainthm}{Main Theorem}
\maketitle
\newcommand{\pp}{\mathfrak{p}}
\renewcommand{\gg}{\mathfrak{g}}
\newcommand{\DD}{\mathcal{D}}
\newcommand{\F}{\ensuremath{\mathbb F}}
\newcommand{\Fp}{\ensuremath{\F_p}}
\newcommand{\Fl}{\ensuremath{\F_l}}
\newcommand{\Fpbar}{\overline{\Fp}}
\newcommand{\Fq}{\ensuremath{\F_q}}
\newcommand{\PP}{\mathbb{P}}
\newcommand{\PPone}{\mathfrak{p}_1}
\newcommand{\PPtwo}{\mathfrak{p}_2}
\newcommand{\PPonebar}{\overline{\PPone}}
\newcommand{\N}{\ensuremath{\mathbb N}}
\newcommand{\Q}{\ensuremath{\mathbb Q}}
\newcommand{\Qbar}{\overline{\Q}}
\newcommand{\R}{\ensuremath{\mathbb R}}
\newcommand{\Z}{\ensuremath{\mathbb Z}}
\newcommand{\SSS}{\ensuremath{\mathcal{S}}}
\newcommand{\Rn}{\ensuremath{\mathbb R^n}}
\newcommand{\Ri}{\ensuremath{\R^\infty}}
\newcommand{\C}{\ensuremath{\mathbb C}}
\newcommand{\Cn}{\ensuremath{\mathbb C^n}}
\newcommand{\Ci}{\ensuremath{\C^\infty}}
\newcommand{\U}{\ensuremath{\mathcal U}}
\newcommand{\gn}{\ensuremath{\gamma^n}}
\newcommand{\ra}{\ensuremath{\rightarrow}}
\newcommand{\fhat}{\ensuremath{\hat{f}}}
\newcommand{\ghat}{\ensuremath{\hat{g}}}
\newcommand{\hhat}{\ensuremath{\hat{h}}}
\newcommand{\covui}{\ensuremath{\{U_i\}}}
\newcommand{\covvi}{\ensuremath{\{V_i\}}}
\newcommand{\covwi}{\ensuremath{\{W_i\}}}
\newcommand{\Gt}{\ensuremath{\tilde{G}}}
\newcommand{\gt}{\ensuremath{\tilde{\gamma}}}
\newcommand{\Gtn}{\ensuremath{\tilde{G_n}}}
\newcommand{\gtn}{\ensuremath{\tilde{\gamma_n}}}
\newcommand{\gnt}{\ensuremath{\gtn}}
\newcommand{\Gnt}{\ensuremath{\Gtn}}
\newcommand{\Cpi}{\ensuremath{\C P^\infty}}
\newcommand{\Cpn}{\ensuremath{\C P^n}}
\newcommand{\lla}{\ensuremath{\longleftarrow}}
\newcommand{\lra}{\ensuremath{\longrightarrow}}
\newcommand{\Rno}{\ensuremath{\Rn_0}}
\newcommand{\dlra}{\ensuremath{\stackrel{\delta}{\lra}}}
\newcommand{\pii}{\ensuremath{\pi^{-1}}}
\newcommand{\la}{\ensuremath{\leftarrow}}
\newcommand{\gonem}{\ensuremath{\gamma_1^m}}
\newcommand{\gtwon}{\ensuremath{\gamma_2^n}}
\newcommand{\omegabar}{\ensuremath{\overline{\omega}}}
\newcommand{\dlim}{\underset{\lra}{\lim}}
\newcommand{\ilim}{\operatorname{\underset{\lla}{\lim}}}
\newcommand{\Hom}{\operatorname{Hom}}
\newcommand{\Ext}{\operatorname{Ext}}
\newcommand{\Part}{\operatorname{Part}}
\newcommand{\Ker}{\operatorname{Ker}}
\newcommand{\im}{\operatorname{im}}
\newcommand{\ord}{\operatorname{ord}}
\newcommand{\unr}{\operatorname{unr}}
\newcommand{\B}{\ensuremath{\mathcal B}}
\newcommand{\Ocr}{\ensuremath{\Omega_*}}
\newcommand{\Rcr}{\ensuremath{\Ocr \otimes \Q}}
\newcommand{\Cptwok}{\ensuremath{\C P^{2k}}}
\newcommand{\CC}{\ensuremath{\mathcal C}}
\newcommand{\gtkp}{\ensuremath{\tilde{\gamma^k_p}}}
\newcommand{\gtkn}{\ensuremath{\tilde{\gamma^k_m}}}
\newcommand{\QQ}{\ensuremath{\mathcal Q}}
\newcommand{\I}{\ensuremath{\mathcal I}}
\newcommand{\sbar}{\ensuremath{\overline{s}}}
\newcommand{\Kn}{\ensuremath{\overline{K_n}^\times}}
\newcommand{\tame}{\operatorname{tame}}
\newcommand{\Qpt}{\ensuremath{\Q_p^{\tame}}}
\newcommand{\Qpu}{\ensuremath{\Q_p^{\unr}}}
\newcommand{\scrT}{\ensuremath{\mathfrak{T}}}
\newcommand{\That}{\ensuremath{\hat{\mathfrak{T}}}}
\newcommand{\Gal}{\operatorname{Gal}}
\newcommand{\Aut}{\operatorname{Aut}}
\newcommand{\tors}{\operatorname{tors}}
\newcommand{\Zhat}{\hat{\Z}}
\newcommand{\linf}{\ensuremath{l_\infty}}
\newcommand{\Lie}{\operatorname{Lie}}
\newcommand{\GL}{\operatorname{GL}}
\newcommand{\End}{\operatorname{End}}
\newcommand{\aone}{\ensuremath{(a_1,\ldots,a_k)}}
\newcommand{\raone}{\ensuremath{r(a_1,\ldots,a_k,N)}}
\newcommand{\rtwoplus}{\ensuremath{\R^{2  +}}}
\newcommand{\rkplus}{\ensuremath{\R^{k +}}}
\newcommand{\length}{\operatorname{length}}
\newcommand{\Vol}{\operatorname{Vol}}
\newcommand{\cross}{\operatorname{cross}}
\newcommand{\GoN}{\Gamma_0(N)}
\newcommand{\GeN}{\Gamma_1(N)}
\newcommand{\GAG}{\Gamma \alpha \Gamma}
\newcommand{\GBG}{\Gamma \beta \Gamma}
\newcommand{\HGD}{H(\Gamma,\Delta)}
\newcommand{\Ga}{\mathbb{G}_a}
\newcommand{\Div}{\operatorname{Div}}
\newcommand{\Divo}{\Div_0}
\newcommand{\Hstar}{\cal{H}^*}
\newcommand{\txon}{\tilde{X}_0(N)}
\newcommand{\sep}{\operatorname{sep}}
\newcommand{\notp}{\not{p}}
\newcommand{\Aonek}{\mathbb{A}^1/k}
\newcommand{\Wa}{W_a/\mathbb{F}_p}
\newcommand{\Spec}{\operatorname{Spec}}

\newcommand{\abcd}{\left[ \begin{array}{cc}
a & b \\
c & d
\end{array} \right]}

\newcommand{\abod}{\left[ \begin{array}{cc}
a & b \\
0 & d
\end{array} \right]}

\newcommand{\unipmatrix}{\left[ \begin{array}{cc}
1 & b \\
0 & 1
\end{array} \right]}

\newcommand{\matrixeoop}{\left[ \begin{array}{cc}
1 & 0 \\
0 & p
\end{array} \right]}

\newcommand{\w}{\omega}
\newcommand{\Qpi}{\ensuremath{\Q(\pi)}}
\newcommand{\Qpin}{\Q(\pi^n)}
\newcommand{\pibar}{\overline{\pi}}
\newcommand{\pbar}{\overline{p}}
\newcommand{\lcm}{\operatorname{lcm}}
\newcommand{\trace}{\operatorname{trace}}
\newcommand{\OKv}{\mathcal{O}_{K_v}}
\newcommand{\Abarv}{\tilde{A}_v}
\newcommand{\kbar}{\overline{k}}
\newcommand{\Kbar}{\overline{K}}
\newcommand{\pl}{\rho_l}
\newcommand{\plt}{\tilde{\pl}}
\newcommand{\plo}{\pl^0}
\newcommand{\Du}{\underline{D}}
\newcommand{\A}{\mathbb{A}}
\newcommand{\D}{\underline{D}}
\newcommand{\op}{\operatorname{op}}
\newcommand{\Glt}{\tilde{G_l}}
\newcommand{\gl}{\mathfrak{g}_l}
\newcommand{\gltwo}{\mathfrak{gl}_2}
\newcommand{\sltwo}{\mathfrak{sl}_2}
\newcommand{\h}{\mathfrak{h}}
\newcommand{\tA}{\tilde{A}}
\newcommand{\sss}{\operatorname{ss}}
\newcommand{\X}{\Chi}
\newcommand{\ecyc}{\epsilon_{\operatorname{cyc}}}
\newcommand{\hatAl}{\hat{A}[l]}
\newcommand{\sA}{\mathcal{A}}
\newcommand{\sAt}{\overline{\sA}}
\newcommand{\OO}{\mathcal{O}}
\newcommand{\OOB}{\OO_B}
\newcommand{\Flbar}{\overline{\F_l}}
\newcommand{\Vbt}{\widetilde{V_B}}
\newcommand{\XX}{\mathcal{X}}
\newcommand{\GbN}{\Gamma_\bullet(N)}
\newcommand{\Gm}{\mathbb{G}_m}
\newcommand{\Pic}{\operatorname{Pic}}
\newcommand{\FPic}{\textbf{Pic}}
\newcommand{\solv}{\operatorname{solv}}
\newcommand{\Hplus}{\mathcal{H}^+}
\newcommand{\Hminus}{\mathcal{H}^-}
\newcommand{\HH}{\mathcal{H}}
\newcommand{\Alb}{\operatorname{Alb}}
\newcommand{\FAlb}{\mathbf{Alb}}
\newcommand{\gk}{\mathfrak{g}_k}
\newcommand{\car}{\operatorname{char}}
\newcommand{\Br}{\operatorname{Br}}
\newcommand{\gK}{\mathfrak{g}_K}
\newcommand{\coker}{\operatorname{coker}}
\newcommand{\red}{\operatorname{red}}
\newcommand{\CAY}{\operatorname{Cay}}
\newcommand{\ns}{\operatorname{ns}}

\begin{abstract}
Fix a non-negative integer $g$ and a positive integer $I$ dividing
$2g-2$.  For any Henselian, discretely valued field $K$ whose
residue field is perfect and admits a degree $I$ cyclic extension,
we construct a curve $C_{/K}$ of genus $g$ and index $I$.  We can
in fact give a complete description of the finite extensions $L/K$
such that $C$ has an $L$-rational point. Applications are
discussed to the corresponding problem over number fields.  S.
Sharif, in his 2006 Berkeley thesis, has independently obtained
similar (but not identical) results.  Our proof, however, is
different: via deformation theory, we reduce to the problem of
finding suitable actions of cyclic groups on finite graphs.
\end{abstract}

\noindent Some terminological conventions: by a variety (resp. a
curve) over a field $K$ we will mean a finite-type $K$-scheme
which is smooth, projective and geometrically integral (resp. of
dimension one). By a variety (resp. a curve) over a field $k$ we
will mean a finite-type $k$-scheme which is geometrically integral
(resp. of dimension one) but possibly incomplete or singular.  If
$V$ is a variety defined over $K$ and $L/K$ is a field extension,
we say that $L$ \textbf{splits} $V$ if $V(L) \neq \emptyset$.

\section{Introduction}
\noindent Given a variety $V$ defined over a field $K$, one would
like to determine whether $V$ has a $K$-rational point, and if it
does not, to say something about $\mathcal{S}(V)$, the set of
finite field extensions $L/K$ for which $V$ acquires an
$L$-rational point.  This is a very difficult problem: e.g., it is
believed by many (but unproved) that there is no algorithm for the
task of
deciding whether a variety $V_{/\Q}$ has a $\Q$-rational point. \\
\indent In order to quantify the second part of the question, it
is natural to introduce the \textbf{index} $I(V)$ of a variety
$V(K)$: it is the greatest common divisor of all degrees of closed
points on $V$ (so $V(K) \neq \emptyset \implies I(V) = 1$, but not
conversely).  For a curve $C_{/K}$, the index is equal to the
least positive degree of a line bundle on $K$.  If $C$ has genus
$g$, then the canonical bundle $\Omega^1_{C/K}$ has degree equal
to $2g-2$ and $I(C) \ | \ 2g-2$.\footnote{Note that this holds --
vacuously -- even when $g =1$.}
\\ \indent
In general, to know $I(V)$ is much less than to know
$\mathcal{S}(V)$: we need not know any particular splitting field
$L/K$, nor even the least possible degree of a splitting field (a
quantity called the \textbf{m-invariant} $m(V)$ in \cite{PLC5}).
For instance, it follows from the Weil bound for curves that every
variety over a finite field has $I(V) = 1$; nevertheless computing
$\mathcal{S}(V)$ or even $m(V)$ is still a nontrivial task.
\\ \\
As usual, when the ``direct problem'' of computation of an
invariant is sufficiently difficult, it is natural to consider as
well the ``inverse problem'' of which invariants actually arise.
In particular, we may well ask:
\begin{ques}
Fix a field $K$.  For which pairs $(g,I) \in \N \times \Z^+$ does
there exist a curve $C_{/K}$ with $I(C) = I$?
\end{ques}
\noindent As above, the existence of the canonical divisor implies
that a necessary condition is $I \ | \ 2g-2$.  This condition is,
of course, not sufficient: e.g., we have $I(C) = 1$ for all curves
when $K$ is finite or PAC,\footnote{A field $K$ is \textbf{P}seudo
\textbf{A}lgebraically \textbf{C}losed if every variety $V_{/K}$
has a rational point. This includes separably closed fields, but
there are many others: see \cite{FJ}.} whereas if $K$ is $\R$ (or
is real-closed, or pseudo-real closed) we have $I(C) \ | \ 2$.
\\ \\
In contemporary arithmetic geometry, the fields of most interest
are those which are infinite and finitely generated (``IFG'').
\begin{conj}
\label{CONJ2} Let $K$ be an IFG field. Then for any $g$ and $I$
with $I \ | \ 2g-2$, there exists a curve $C_{/K}$ of genus $g$
with $I(C) = I$.
\end{conj}
\noindent Remark 1.1: Conjecture \ref{CONJ2} is true when $g = 0$;
this is a relatively easy exercise involving quaternion algebras
which we leave to the reader.
\\ \\
Remark 1.2: The main result of \cite{PLC3} is that Conjecture
\ref{CONJ2} holds for $g=1$ when $K$ is a number field.  It ought
to be possible to extend this result to positive characteristic
global fields (using the period-index obstruction map in flat
cohomology) and then to all IFG fields (using isotrivial elliptic
curves). Indeed, in genus one it is natural to make a much
stronger conjecture: see \cite[Conjecture 1]{Clark-Sharif}.
\\ \\
Let us now present evidence for Conjecture \ref{CONJ2} for curves
of higher genus ($g \geq 2$).  The following result uses the
author's work on the genus one case together with some simple
covering considerations (essentially those suggested to the author
by Bjorn Poonen in 2003) to attain a solution for ``small
indices.''
\begin{thm}
\label{PTHM} Let $K$ be a number field, $g \in \N$ and $k \in
\Z^+$ with $k \ | \ g-1$.  Then there exists a curve $Y_{/K}$ of
genus $g$ and index $I = \frac{g-1}{k}$.
\end{thm}
\noindent Proof: Since $K$ has characteristic different from $2$,
it is especially easy to see that there exist (hyperelliptic)
curves over $K$ of all genera with $K$-rational (Weierstrass)
points; we may therefore assume that $I > 1$.  By Remark 1.2 there
exists a curve $X_{/K}$ of genus one and index $I$. By \cite[Prop.
11]{PLC7}, there exist linearly equivalent $K$-rational divisors
$D_1$ and $D_2$ on $C$, both effective of degree $kI$, such that
the support of $D_1-D_2$ has cardinality $2kI$.  Let $f \in K(X)$
be the rational function (uniquely determined up to a constant
multiple) with divisor $D_1 - D_2$, and let $\varphi: Y \ra X$ be
the branched covering corresponding to the extension of function
fields $K(X)(\sqrt{f})/K(X)$.  By the Riemann-Hurwitz formula, $Y$
has genus $g =kI + 1$.  As for any covering of curves, we have
$I(Y) \geq I(X) = I$; conversely, any point $P$ in the support of
$D_1$ has degree $I$ and is a ramification point for $\varphi$, so
its unique preimage $\tilde{P} \in Y(\overline{K})$ has degree
$I$. Thus we have $I(Y) = m(Y) = I$.
\\ \\
In this note we are mainly concerned with presenting a technique
for constructing curves $C$ (and perhaps higher-dimensional
varieties) over \emph{local} fields for which the set
$\mathcal{S}(C)$ can be explicitly computed. Here is our main
result:
\begin{mainthm}
\label{MT} Let $K$ be a Henselian discretely valued field with
perfect residue field.  Assume that there is a cyclic, degree $I$,
unramified extension $K_I/K$ (``hypothesis $K_I$''). Then for any
non-negative integer $g$ such that $I \ | \ 2g-2$, there
exists a curve $C_{/K}$ with the following properties: \\
a) If $I$ is odd or $g = 1$, then a finite extension $L/K$ splits
$C$ iff $L
\supset K_I$. \\
b) Otherwise, a finite extension $L/K$ splits $C$ iff either of
the following holds: \\
$\bullet$ $L \supset K_I$; or \\
$\bullet$ $L \supset K_{I/2}$ -- the unique subextension of
$K_{I}/K$ of degree $I/2$ -- and $2 \ | \ e(L/K)$.
\end{mainthm}
\noindent Remark 1.3: So if $K$ satisfies hypothesis $K_I$ for all
$I$, then for all $I \ | \ 2g-2$ there is a curve $C_{/K}$ of
genus $g$ and index $I$.  In particular this holds when the
residue field is \emph{finitely generated}, when we may construct
$K_I$ by adjoining suitable roots of unity.
\\ \\
Remark 1.4: Some hypothesis on the existence of unramified
extensions is necessary, since if $k$ is algebraically closed,
then $I(C) \ | \ g-1$ \cite[Remark 1.8]{BL}. On the other hand,
one could ask for a classification of all possible indices under
the milder assumption that $K$ admits a quadratic unramified
extension, e.g. in the case $K = \R((T))$. We shall not pursue
this here.
\begin{cor}
\label{COR} Let $K$ be an infinite, finitely generated field.  For
any $g \in \N$, there exists a finite extension $L/K$ and a genus
$g$ curve $C_{/L}$ with $I(C) = 2g-2$.
\end{cor}
\noindent \emph{Proof} of Corollary \ref{COR}: Such a field $K$
admits a discrete valuation $v$ with finitely generated residue
field $k$, so by Remark 1.3 we may apply our Main Theorem to the
Henselization $K_v$ of $K$.  We thus get a curve $C$ of genus $g$
and index $I$ defined over $K_v$, which is an algebraic extension
of $K$; it is evidently defined over some finite extension $L$ of
$K$.
\\ \\
Remark 1.5: Only a few days after the results of this paper were
first obtained I received a copy of the 2006 Berkeley thesis of
S.I. Sharif \cite{Sharif}, which contains closely related results.
\begin{thm}(Sharif,  \cite{Sharif})
\label{SHARIFTHM}\\
a) Let $K$ be a locally compact discretely valued field of
characteristic different from $2$.  Then for any $(g,I) \in \N
\times \Z^+$ with $I \ | \ 2g-2$, there exists a curve $C_{/K}$ of
genus $g$ and index $I$. \\
b) For any $g \in \N$ and $I \in \Z^+$ such that $4 \not{|} \ | \
I \ | \ 2g-2$, there exists a number field $K = K(g,I)$ and a
curve $C_{/K}$ of genus $g$ and index $I$.
\end{thm}
\noindent By Remark 1.3, part a) of Theorem \ref{SHARIFTHM} is a
special case of our Main Theorem, whereas part b) is similar in
spirit to Theorem \ref{PTHM} and Corollary \ref{COR} but not
directly comparable to either one.  On the other hand, Sharif's
results go further than those presented here in that he also
considers the possible values of the \emph{period} $P$ (the least
positive degree of a $K$-rational divisor class), and -- comparing
with the restrictions on period and index obtained by Lichtenbaum
-- his constructions give the complete list of possible values of
$(g,P,I)$ for curves over
locally compact fields of characteristic different from $2$.  \\
\indent The strategy of Sharif's proof -- namely, construction of
degree two covers of curves of genus one and two via $p$-adic
theta functions -- is quite different from the proof of our Main
Theorem.
\\ \\
Let us now say something about how we shall prove the Main
Theorem.  There are three steps: (i) we observe that if a curve
$C_{/K}$ admits a regular model $C_{/R}$ with \emph{semistable}
special fiber, then the set $\mathcal{S}(C)$ of finite extensions
$L/K$ with $C(L) \neq \emptyset$ depends only on the special
fiber. This reduces us to the problem of determining rational
points on a semistable curve over the special fiber, and the next
idea is that the ``more singular'' the special fiber, the
\emph{easier} it is to analyze its set of rational points.  In
particular, if we restrict to the case of \emph{totally
degenerate} special fibers -- in which each geometric component
has genus zero -- then the information we want can be gleaned from
the action of the Galois group $\mathfrak{g}_k =
\Gal(\overline{k}/k)$ on the \emph{dual graph} whose vertices are
the components and edges are the intersection points.  (ii) A
fundamental result in the deformation theory of one-dimensional
complete intersections (due essentially to Grothendieck) tells us
that \emph{every} semistable curve over $k$ lifts to a regular
scheme over the valuation ring $R$ with smooth generic fiber;
complementing this with the elementary observation (due to A.
P\'al) that every cubic graph is the dual graph of a semistable
curve over $k$, we are reduced to (iii) a combinatorial problem
involving the construction of a family of finite graphs with
suitable Euler characteristic and automorphisms by a finite cyclic
group.  It turns out that one solution to the problem is obtained
by using a well-known family of graphs, the \emph{M\"obius
ladders}.  It seems enlightening to present the construction in
terms of Cayley graphs, and we do so here, although we have
written up the endgame in such a way as to offer the reader
unfamiliar with the formalism of Cayley graphs a choice: he may
either learn this materal here or bypass it in favor of a simple,
concrete description.
\\ \\
It is natural to try to generalize the results to a broader class
of algebraic varieties, e.g. varieties $V_{/K}$ admitting a
regular SNC model (in particular all algebraic curves), or
varieties whose special fiber is completely degenerate.  I also
suspect that there should be a more elegant approach via
non-Archimedean uniformization.
\\ \indent
I have decided not to try to work out these generalizations in the
present note, because (i) it seems they would require
significantly more technical apparatus, whereas our present
methods -- with the exception of the aforementioned
deformation-theoretic result that we treat as a ``black box'' --
are rather elementary; and (ii) because our work overlaps
substantially with S. Sharif's thesis, it seems best to record
once and for all our work that was done independently of
\cite{Sharif}, so that our subsequent work can draw freely on both
sources.

\section{Local points on semistable curves}
\noindent The results of this section may be well-known to some,
but -- especially in the absence of any satisfactory reference --
we prefer to work them out in detail.
\\ \\
Let $C_{/K}$ be a curve.  Suppose that $C$ has semistable
reduction: that is, there exists a regular arithmetic surface $C$
over the valuation ring $R$ of $K$ with generic fiber isomorphic
to $C$ and with special fiber a semistable curve $C_{/k}$.  Recall
that a semistable curve $C_{/k}$ is a one-dimensional projective
$k$-scheme which is reduced, geometrically connected, and whose
only singularities are ordinary double points.
  \\ \indent Write $C_{/\overline{k}} =
\sum_{i=1}^N C_i$, so that the $C_i$ are the
$\overline{k}$-irreducible components (which we will henceforth
call the components).
\\ \indent
There is a natural $\gk$-action on the set of components.  We will
say that a component $C_i$ is \emph{defined} over a finite field
extension $l/k$ if $\mathfrak{g}_l$ fixes $C_i$.  There is
evidently a unique minimal such field extension (necessarily
Galois over $k$), which we denote by $l_i$. Let $d_i = [l_i:k_i]$
and $d = \gcd d_i$.
\\ \\
For any reduced finite-type scheme $S_{/k}$, define its
\textbf{nonsingular index} $I^{\ns}(S)$ to be the index of the
nonsingular locus $S^{\ns}$.  Put $I_i := I^{\ns}(C_i)$.
\begin{thm}
\label{INDEXTHM} Let $C_{/R}$ be a regular arithmetic surface with
generic fiber $C_{/K}$ and semistable special fiber $C_{/k}$. Then
\[I(C_{/K}) = I^{\ns}(C_{/k}) = \gcd_i (d_i \cdot I_i). \]
\end{thm}
\noindent The proof will come later in this section.
More information can be obtained under the following hypothesis:
\\ \\
(A) For every finite extension $l/k$, every component $C_i$ which
is defined over $l$ has an $l$-rational point which is not a nodal
point of $C$.
\begin{prop}
\label{INDEXPROP}
Maintain the hypotheses of Theorem
\ref{INDEXTHM}, and assume also
(A).  \\
a) If $d$ is odd, then a finite extension $L/K$ splits $C$ iff its
maximal unramified extension $L'/K$ splits $C$ iff its residue extension contains $l_i$ for some $i$. \\
b) If $d$ is even, then a finite extension $L/K$ splits $C$ iff
either \\
(i) the residue extension $l$ contains $l_i$ for some
$i$; or  \\
(ii) the residue extension $l/k$ is such that $\mathfrak{g}_{l}$
stabilizes a pair of intersecting components, and $e(L/K)$ is
even.
\end{prop}
\noindent Remark: Let us say that we are in \textbf{Case 1} if $d$
is odd or condition (bii) does not occur, and that we are in
\textbf{Case 2} if $d$ is even and condition (bii) occurs.
\\ \\
It will be convenient to introduce the (so-called) dual graph
$\mathcal{G} = (\mathcal{V},\mathcal{E})$, an undirected,
connected, finite graph whose vertices are the components of $C$,
and where vertices $C_i$ and $C_i$ are linked by $C_i \cdot C_j$
edges.  The natural action of the Galois group $\gk$ on components
and on singular geometric points gives rise to an action of $\gk$
on $\mathcal{G}$ by graph-theoretical automorphisms.
\\ \\
\emph{Proof} of Proposition \ref{INDEXPROP}: By a well-known
version of Hensel's Lemma (e.g. \cite{JL}), $C(K) \neq \emptyset$
iff $C$ has a smooth $k$-rational point.  In the present case this
occurs iff some vertex $C_i$ of $\mathcal{G}$ is fixed by the
action of $\gk$. Since regularity of a model is unaffected by
unramified base change, an unramified extension $L/K$ splits $C$
iff its residue extension $l$ contains $l_i$ for some $i$.
\\ \indent
Let us now consider the effect of making a totally ramified base
extension $L/K$.  We can form the arithmetic surface $C'_{/S} := C
\otimes_R S$ (where $S$ is the valuation ring of $R$), and the
special fiber is still $C_{/l} = C_{/k}$, but $C'$ will no longer
be regular (unless $C_{/k}$ was smooth), because the complete
strict local ring at a singularity will now be of the form
$S^{\unr}[[x,y]]/(xy-\pi^e)$, where $\pi$ is a uniformizer for $L$
and $e = e(L/K)$ is the relative ramification index.  The remedy
is well-known -- we must blowup $e-1$ times to replace the
intersection point with a chain of $e-1$ rational curves. However,
keep in mind that we are really blowing up a closed point whose
residue field may be larger than $k$, or in other words, we are
simultaneously blowing up each point in the $\gk$-orbit of the
given singular point, and there is, in an evident way, an induced
$\gk$-action on these chains.  In order for there to be a smooth
$l$-rational point \emph{after} the blowing-up process which was
not there \emph{before} that process, necessary and sufficient
conditions are: first, for some chain to be $\gk$-stable -- in
other words we need for a pair $\{C_i, \ C_j\}$ of intersecting
components to be preserved by $\gk$; and \\
second, for the chain to have odd length, i.e., for $e$ to be
even.  This completes the proof of Proposition \ref{INDEXPROP}.
\\ \\
Before beginning the proof of Theorem 6 we will record two quick
lemmas.
\begin{lemma}
\label{LEMMA1} If $V_{/k}$ is a finite-type reduced scheme and
$l/k$ is a finite field extension, then $I(V_{/k}) \ | \ [l:k]
 \cdot I(V_{/l})$.
\end{lemma}
\noindent Proof: Let $D_l$ be an $l$-rational zero-cycle on $V$ of
degree $I(V_{/l})$; its trace from $l$ down to $k$ is a
$k$-rational zero-cycle of degree $[l:k] \cdot I(V_{/l})$.
\begin{lemma}
\label{LEMMA2} The nonsingular index $I^{\ns}(V)$ of a variety
$V_{/k}$ is a birational invariant.  In particular the nonsingular
index $I^{\ns}(C)$ of a curve $C_{/k}$ is unchanged by the removal
of finitely many closed points.
\end{lemma}
\noindent Proof: If $k$ is finite, then much more is true: it
follows from the Weil bounds for curves over finite fields that
for all finite-type geometrically integral schemes $V_{/\F_q}$
$I^{\ns}(V) = 1$.  The argument is well-known to field
arithmeticians: see \cite{FJ}.  When $k$ is infinite, see
\cite[p.8]{CT}.
\\ \\
Remark 2.1: With $I$ instead of $I^{\ns}$, the conclusion of Lemma
\ref{LEMMA2} does not follow: take $X^2+Y^2+Z^2 = 0$ over $\R$.
\\ \\
We come now to the proof of Theorem \ref{INDEXTHM}: \\ \\  Step 1:
We will show that $I^{\ns}(C_{/k}) = \gcd_i d_i I_I$.
 Indeed for any $i$, $1 \leq i \leq N$
\[I^{\ns}(C_{/k}) = I(C^{\ns}_{/k}) \ | \ [l_i:k] \cdot
I(C^{\ns}_{/l_i}) \ | \  d_i I_i \] (these divisibilities use
Lemmas \ref{LEMMA1} and \ref{LEMMA2}), so $I^{\ns}(C_{/k}) \ | \
\gcd d_i I_i$. For the converse, $I^{\ns}(C_{/k})$ is the gcd of
all degrees of field extensions $l/k$ such that $C$ has a smooth
$l$-rational point $P$.  Thus $l$ must be a field of definition
for the component on which $P$ lies -- i.e., $l_i \subset l$, and
then clearly $I_i \ | \ [l:l_i]$.
\\ \\
Step 2: It is clear from Hensel's Lemma that $I(C_{/K}) \ | \
I^{\ns}(C_{/k})$; more precisely, the fields $l$ for which $C$
acquires a smooth $l$-rational point correspond to the unramified
splitting fields.  It remains to account for the possibility that
$L$ splits $K$ when its maximal unramified subextension $K'$ does
not.  As in the proof of Proposition \ref{INDEXPROP}, this can
only happen when there is a $k'$-rational singular point $q$  on
$C$ and $e(L/K)$ is even.  Suppose first that $q$ is the
intersection of distinct components $\{C_i, \ C_j\}$.  Then the
pair of components is stabilized by $\mathfrak{g}_{k'}$, and $q$
is a smooth point on each component, so $I^{\ns}(C_{/k'}) \ | \
2$, so that there is an unramified splitting field $K''/K$ with
$[K'':K] \ | \ [L:K]$.  The other possibility is that $q$ is a
nodal singularity on a single component $C_i$ which is defined
over $k$.  But then the preimage of $q$ in the normalization
$\tilde{C_i}$ of $C_i$ is a rational divisor of degree $2$, which
by Lemma \ref{LEMMA2} can be ``moved'' to a rational divisor of
degree $2$ with support disjoint from the singular locus, hence
projecting down to a rational divisor of degree $2$ on $C_i$,
i.e., $I_i \ | \ 2$ and again we have found an unramified
splitting field $K''$ such that $[K'':K] \ | \ [L:K]$.  This
completes the proof of Theorem \ref{INDEXTHM}.
\\ \\
Remark 2.2: The proof of Theorem \ref{INDEXTHM} shows that curves
$C_{/K}$ with semistable reduction are ``unramified'' in a strong
sense: not only do they possess unramified splitting fields, but
also the index can be calculated using only unramified extensions.
However, in the absence of Hypothesis (A), the \emph{least} degree
of a field extension $L/K$ such that a curve $C_{/K}$ with
semistable reduction acquires rational points will not in general
be attainable by an unramified extension.  See \cite{PLC5} for a
naturally occurring example (involving the Shimura curves
$X^D_0(N)$) where in order to show that a semistable curve has
points over a quartic base extension, a ramified base change was
used in a perhaps essential way.
\section{Proof of the Main Theorem}
\noindent As the reader has probably guessed, the curves $C_{/K}$
referred to in the Main Theorem will be such that they admit
regular models whose special fibers are semistable and satisfy
hypothesis $(A)$, so information about their splitting fields will
reduce to an analysis of the dual graph.  In fact, we we will
place ourselves in a situation in which we need only construct the
dual graph and not the arithmetic surface itself.  This is done
via the following two results:
\begin{thm}
\label{DEFTHM}
For any semistable curve $C_{/k}$, there exists a
regular arithmetic surface whose special fiber is isomorphic to
$C$ and whose generic fiber is a curve over $K$.
\end{thm}
\noindent \emph{Proof}: This is (I gather) a standard result in
deformation theory.  A relatively accessible reference is
\cite[4.4]{Vistoli}.
\\ \\
Thus it suffices to construct singular curves over the residue
field $k$.  We will in fact construct \emph{totally degenerate}
semistable curves, namely with each component of geometric genus
$0$.  For this:
\begin{lemma}
\label{LEMMA10} Let $\mathcal{G}$ be any connected graph in which
each vertex has degree at most $3$.  Let $G$ be a finite group
acting on $\mathcal{G}$ by automorphisms.  Given a field $k$, a
Galois extension $l/k$ and an isomorphism $\mathfrak{g}_{l/k}$
with $G$, there is a totally degenerate semistable curve $C_{/k}$
whose dual graph is isomorphic to $\mathcal{G}$, under an
isomorphism which identifies the Galois action on $\mathcal{G}$
with the action of $G$.
\end{lemma}
\noindent Proof: This is shown in \cite{Pal} under the hypothesis
that $k$ is infinite but without the hypothesis that $\mathcal{G}$
have degree $3$.  The infinitude of $k$ is used precisely to
ensure that the intersection points of the graph can be identified
with $k$-rational points of $\PP^1(k)$.  Since $\# \PP^1(k) \geq
3$ for all $k$, the argument goes through verbatim with the
hypothesis of degree at most $3$.
\\ \\
Recall that the arithmetic genus of a totally degenerate
semistable curve $C_{k}$ (which is the genus of any smooth lift
$C_{/K}$ in the usual sense) is just $1 - \chi(\mathcal{G})$,
where $\chi$ is the Euler characteristic of the dual graph in the
usual topological sense, computable as the number of vertices
minus the number of edges.
\\ \\
The curves constructed by Lemma \ref{LEMMA10} satisfy hypothesis
(A) unless the residue field $k$ is $\F_2$.  More precisely, what
we need is that for each finite extension $l/k$, every component
which is defined over $l$ has at most $\# l$ singular points.
Since our graphs have degree at most $3$, the only problematic
case is when $k = \F_2$ and $I = 1$ (because if $I > 1$, we only
want points over an extension with larger residue field).  But
this is a trivial case: it is enough, for instance, to find a
nonsingular curve $C_{/\F_2}$, of genus $g$, and with $C(\F_2)
\neq \emptyset$.  Or, staying with the same graph-theoretical
strategy, we need only to find, for all $g \geq 0$, a connected
graph with Euler characteristic $1-g$, in which each vertex has
degree at most $3$, and at least one vertex has degree at most
$2$.  Of course such graphs exist: for $g = 0$ take the graph with
one vertex and no edges (the dual graph of $\PP^1$), and for $g
\geq 1$ we can build such a graph out of $g$ ``coathangers'' (the
graph with vertex set $\{0,1,2,3\}$ and $0 \sim 1, \ 0 \sim 2, \ 0
\sim 3, \ 2 \sim 3$).  Henceforth we will assume that $I > 1.$
\\ \\
Let $G$ be a group and $S \subset G$ such that $S = S^{-1}$, $1
\not \in S$,  and $\langle S \rangle = G$.  We define the Cayley
graph $\CAY(G,S)$, a simple (no loops, no multiple edges)
undirected graph whose vertex set is $G$ itself, and with
\[g \sim g' \iff \exists s \in S \ | \ gs = g'.
\]
Note the following (almost tautological) properties of
$\CAY(G,S)$: (i) it is connected; each vertex has degree $\# S$;
(iii) it admits a left $G$-action which is free on vertices, and
free on edges unless $S$ contains an element of order $2$.\\
(iv) If $G$ is finite,
\[\chi(\CAY(G,S)) =  \#G \left(1-\frac{\#S}{2}\right). \]
(v) If $\rho: H \hookrightarrow G$ is an embedding, then $H$ acts
on $\CAY(G,S)$, freely on vertices and freely on edges unless
$\rho(H) \cap S$ contains an element of order $2$.
\\ \\
Now let $G_I = \langle \sigma \ | \ \sigma^I = 1 \rangle$, and
identify $G_I$ with $\mathfrak{g}_{K_I/K} = \mathfrak{g}_{k_I/k}$.
\\ \\
When $g = 0$ and $I = 2$, we can take $\mathcal{G} =
\CAY(G_2,\{\sigma\})$, the unique connected graph with two
vertices and one edge: $\chi = 1$.  Since the generator has order
$2$, the unoriented edge gets stabilized, so by Proposition 4 we
get the ``Case 2'' splitting behavior indicated in the theorem.
\\ \\
When $g = 1$ and $I > 2$, we take $\mathcal{G} = \CAY(G,\{\sigma,
\sigma^{-1}\})$, the $I$-cycle: $\chi = 0$.  Here $G$ acts freely
on the edges, so by Proposition 4 we get the ``Case 1'' splitting
behavior indicated in the theorem.
\\ \\
When $g = 1$ and $I = 2$ we can take $\mathcal{G}$ to be the
$2$-cycle, and let $G_I$ act by ``180 degree rotation'', i.e., by
swapping both vertices and both edges: $\chi = 0$, Case 1.  This
graph is nonsimple and \emph{a fortiori} not a Cayley graph
according to our setup.\footnote{On the other hand, the Cayley
graph construction admits several variants.  The $2$-cycle is a
Cayley graph for $G_2$ according to the conventions of e.g.
\cite{delaHarpe}.}  If we insist on seeing a Cayley graph
construction, fix $N
> 1$ and let $\rho: G_2 \hookrightarrow G_{2N}$ be the embedding
$\sigma \mapsto \sigma^N$.
\\ \\
When $g > 1$ and $I = 2g-2$, take $\mathcal{G} =
\CAY(G_I,\{\sigma,\sigma^{-1}, \sigma^{g-1} \}).$  Or, in plainer
terms, start with the $2g-2$ cycle and connect each pair of
antipodal points by an edge: these ``spokes'' do not ruin the
obvious $G_{2g-2}$ action by rotations.  \\ \indent This second
description is graph-theoretically correct (which is, of course,
all that matters for us) but geometrically wrong: the graph does
not really live in the Euclidean plane because the spokes would
have to meet at the center of the circle, adding an unwanted (and
$G_I$-fixed) vertex.  Indeed, when $g = 4$ the graph is precisely
the complete bipartite graph $K_{3,3}$, and for larger $g$ the
graph contains $K_{3,3}$ as a topological subgraph, so these
graphs are \emph{not} embeddable in the Euclidean
plane!\footnote{A real algebraic geometer might be tempted to
point out that the graph naturally lives in the blowup of $\R
\PP^2$ at a single point.}  So here is a ``better'' geometric
description: take a rectangle of length $g$ and height $1$,
subdivide the top and bottom sides into $g-1$ equal parts, and
draw in the $g+1$ equidistant vertical lines linking each vertex
on the top to its corresponding bottom vertex. The resulting graph
has $2g$ vertices and $g+2(g-1)$ edges.  Now identify the right
and left sides of the rectangle with a half-twist, getting a graph
embedded isometrically into the M\"obius band with $2g-2$ vertices
and $g+2(g-1) -1 = 3g-3$ edges with a natural action of the cyclic
group $G_I$ by unit length horizontal rotations. In any case we
have $\chi = 1-g$, and $G_{I}$ acts freely on vertices but not on
edges, Case 2.
\\ \\
When $g > 1$ and $1 < I \ | \ 2g-2$, take the above graph and the
embedding \\ $\rho: G_I \hookrightarrow G_{2g-2}$, $\sigma \mapsto
\sigma^{\frac{2g-2}{I}}$.  If $I$ is odd, we are in Case 1; if $I$
is even, Case 2.
\section{Final remarks}
\noindent Remark 4.1: The genus $0$ case can also be handled using
the correspondence between genus $0$ curves and quaternion
algebras together with the structure theorem for Brauer groups
over Henselian fields with perfect residue field \cite[Ch.
XIII]{CL}.  It is somewhat amusing to remark that in the case of a
finite residue field, our analysis gives a geometric proof of the
well-known fact that a quaternion algebra over a locally compact
field is split by every quadratic extension.  It seems interesting
that the behavior in genus zero, which is often seen as an
anomalous case -- e.g., by virtue of the nonuniqueness of a
minimal regular model -- is completely in line with the behavior
for $g \geq 2$. It is the case of genus one which is truly
exceptional.
\\ \\
Remark 4.2: Consider now the case of a curve $C_{/K}$ admitting an
$R$-model such that the \emph{reduced subscheme} of the special
fiber is semistable: such an arithmetic surface $C_{/R}$ is called
an \textbf{SNC model} of its generic fiber $C_{/K}$.  (It is known
that every curve admits a regular SNC model.)  In this case,
writing $e_i$ for the multiplicity of $C_i$, I find it quite
likely that Theorem 2 generalizes in the form
\[I(C_{/K}) = \gcd_i (d_i \cdot e_i \cdot I_i).\]  At least in the
case of finite residue field (so $I_i = 1$ for all $i$), this
formula follows by combining work of \cite{CTS} and \cite{BL} (as
is observed in \cite[p. 22]{PS}).  As far as I can see, the proofs
-- which are more technically elaborate than in the semistable
case -- should go through for any perfect residue field. However,
unlike the case of semistable reduction, the special fiber does
not itself determine the set $\mathcal{S}(C)$ of all possible
splitting fields.
\\ \\
Remark 4.3: The proverbial alert reader will have noticed that we
defined the index $I(V)$ only for a nonsingular, projective,
geometrically irreducible variety $V_{/K}$.  The definition that
we gave makes sense for any finite-type $K$-scheme $S$ but -- in
view of Remark 4.2 and the fact that this quantity depends only on
the reduced subscheme $S^{\red}$ -- seems to be the wrong
definition for singular schemes.  From the arguments of $\S 2$ --
and in particular the proof of Theorem \ref{INDEXTHM} -- it
follows easily that if $S_{/K}$ is the generic fiber of a regular
model $S_{/R}$ whose special fiber $S_k = \bigcup_i S_i$ is
semistable (i.e., is a reduced finite union of varieties whose
singular points are analytically isomorphic to transversely
intersecting hyperplanes), then $I(S_{/K}) = \gcd_i(d_i I_i)$
where $S_i$ is defined over $l_i$, $d_i = [l_i:k_i]$, and $I_i$ is
the index of $(S_i)_{/l_i}$.  In other words, $\gcd_i(d_i I_i)$ is
a suitable definition of the index of a reduced finite-type scheme
over a perfect field.
\\ \\
As mentioned above in the case of curves, it seems quite plausible
that if $S_{/k}$ is a SNC $k$-scheme such that $S_{/k} = \bigcup_i
S_i$ and $e_i$ is the multiplicity of the component $S_i$, then
defining the index $I'(S_{/k})$ to be $\gcd_i (d_ie_i I_i)$, then
$I(S_{/K}) = I'(S_{/k})$, i.e., the index is preserved by
specialization when the model is regular, the generic fiber is
smooth and the special fiber is SNC.  One wonders about a more
general ``index specialization theorem.''
\begin{conj}
Assume $K$ is a Henselian discrete valuation field with perfect
residue field $k$.  Let $V_{/K}$ be a smooth, projective,
geometrically integral finite-type $K$-scheme and assume that $V$
admits some regular $R$-model, with special fiber $V_{/k}$. Then
$I(V_{/K}) = I'(V_{/k})$.
\end{conj}
\noindent Remark 4.4: Suppose that $V_{/K}$ is an $n$-dimensional
variety admitting a regular model whose special fiber is
completely degenerate, in the sense that it is reduced and every
geometric component is isomorphic to $\PP^n$, and such that the
$k$-fold intersections are transverse and isomorphic to a disjoint
union of $\PP^{n-(k-1)}$'s.  Such varieties arise naturally: they
include all varieties uniformized by Drinfeld's $n$-dimensional
upper halfspace.  For such a variety one can compute the set
$\mathcal{S}(V)$ in terms of the combinatorics of the dual
(simplicial) complex, a generalization of the dual graph of a
curve.  What we lack in this context is an analogue of Theorem
\ref{DEFTHM}: in general, higher-dimensional varieties (even
smooth ones) do not lift smoothly to characteristic $0$. How to
choose the combinatorial geometry of the special fiber so as to
permit a smooth lifting as well as some description of which
classes of higher-dimensional varieties can have completely
degenerate reduction are interesting questions to which we hope to
return in a later work.
\\ \\
\noindent Acknowledgements: This paper records work done during
the 2006 Rational and Integral Points program at MSRI. I thank
MSRI for its hospitality and many of the participants for helpful
conversations, especially Jean-Louis Colliot-Th\'el\`ene, David
Harbater, Bjorn Poonen and Olivier Wittenberg.  Thanks also to
Dino Lorenzini and Ravi Vakil for helping me chase down the
necessary deformation-theoretic facts and to S.I. Sharif for
providing me with an advance copy of his thesis.

\
\end{document}